\documentclass[11pt]{article}
\usepackage{authblk}
\usepackage{graphicx}
\usepackage{color}
\usepackage{amsmath}
\usepackage{amssymb}
\usepackage{amscd}
\usepackage{framed}
\usepackage{bbm}
\usepackage{tcolorbox}
\usepackage{fancyhdr,a4wide}
\usepackage{amsthm}

\newcommand{\R}{\mathbb{R}}
\newcommand{\inr}[1]{\left\langle #1 \right\rangle}

\newcommand{\N}{\mathbb{N}}

\newcommand{\E}{\mathbb{E}}

\newcommand{\PP}{\mathbb{P}}

\newcommand{\eps}{\varepsilon}

\newcommand{\sign}{\operatorname{sign}}

\newtheorem{Theorem}{Theorem}[section]
\newtheorem{Lemma}[Theorem]{Lemma}

\newtheorem{Definition}[Theorem]{Definition}

\newtheorem{Corollary}[Theorem]{Corollary}
\newtheorem{Remark}[Theorem]{Remark}

\numberwithin{equation}{section}

\def \proof {\noindent {\bf Proof.}\ \ }

\def \endproof
{{\mbox{}\nolinebreak\hfill\rule{2mm}{2mm}\par\medbreak}}
\def\IND{\mathbbm{1}}

\newcommand{\bP}{\mathbb{P}}

\def\IND{\mathbbm{1}}

\title{Sharp estimates on random hyperplane tessellations}

\date{}

\author{Sjoerd Dirksen\footnote{Utrecht University, Mathematical Institute (s.dirksen@uu.nl)}, Shahar Mendelson\footnote{University of Warwick, Department of Statistics and The Australian National University, Centre for Mathematics and its Applications (shahar.mendelson@gmail.com)}, Alexander Stollenwerk\footnote{UCLouvain, ICTEAM Institute (alexander.stollenwerk@uclouvain.be)} }

\begin{document}

\maketitle

\begin{abstract}
\noindent We study the problem of generating a hyperplane tessellation of an arbitrary set $T$ in $\R^n$, ensuring that the Euclidean distance between any two points corresponds to the fraction of hyperplanes separating them up to a pre-specified error $\delta$. We focus on random gaussian tessellations with uniformly distributed shifts and derive sharp bounds on the number of hyperplanes $m$ that are required. Surprisingly, our lower estimates falsify the conjecture that $m\sim \ell_*^2(T)/\delta^2$, where $\ell_*^2(T)$ is the gaussian width of $T$, is optimal.
\end{abstract}

\section{Introduction}

We consider the problem of uniformly tessellating sets in a high-dimensional Euclidean space using random hyperplanes. In what follows we use the following notation. For any $a\in \R^n$ and $\tau\in \R$ let
$$
H(a,\tau) = \{x\in \R^n \ : \ \langle a,x\rangle+\tau=0\}
$$
denote the hyperplane defined by $(a,\tau)$. Each hyperplane $H(a,\tau)$ divides $\R^n$ into two parts---the positive and the negative ones, and as a result a collection of hyperplanes endows a tessellation of any given subset $T\subset \R^n$.
 
Let $A:\R^n\to \R^m$ be the matrix with rows $a_1^T,\ldots,a_m^T$, let $\tau=(\tau_1,\ldots,\tau_m)\in \R^m$, and let $\sign:\R^m\to \{-1,1\}^m$ be defined by $\sign(y)=(\sign(y_i))_{i=1}^m$. In this notation, for any given $x\in T$, the bit string $\sign(Ax+\tau)\in \{-1,1\}^m$ encodes the ``cell" in which $x$ is located. Moreover, the normalized Hamming distance
$$
\frac{1}{m}d_H(\sign(Ax+\tau),\sign(Ay+\tau))
$$
is the fraction of the $m$ hyperplanes that separate $x$ and $y$.
\begin{Definition} \label{def:tess}
The hyperplanes $H(a_i,\tau_i)$, $i=1,\ldots,m$, endow a $\delta$-tessellation of the set $T$ if, for some constant $\kappa$,
\begin{equation}
\label{eqn:deltaTessProp}
\left|\frac{\kappa}{m}d_H(\sign(Ax+\tau),\sign(Ay+\tau)) - \|x-y\|_2\right|\leq \delta, \qquad \text{for every} \  x,y \in T.
\end{equation}
\vspace{0.1cm}
\end{Definition}

In other words, a $\delta$-tessellation implies that the Euclidean distance between any two points is proportional to the fraction of hyperplanes separating the two, up to a small additive error of $\delta$. In particular, the Euclidean diameter of all the cells of a $\delta$-tessellation of $T$ is at most $\delta$. Phrased differently, such a $\delta$-tessellation means that the map $f:\R^n\to \{-1,1\}^m$, $f(t)=\sign(At+\tau)$, is a \emph{$\delta$-embedding}; an almost-isometric embedding of $T$ into the discrete cube, with the latter endowed with the metric $\frac{\kappa}{m}d_H$. It is important to stress that an `almost-isometric' embedding here means isometric up to an \emph{additive error}, rather that a multiplicative one.
\par
Apart from being of intrinsic geometric interest, random hyperplanes tessellations and their associated random embeddings appear naturally in at least two applications. First, a $\delta$-embedding can clearly be used to compress a dataset $T$ into a small number of bits. Property \eqref{eqn:deltaTessProp} ensures that one can quickly retrieve information on Euclidean distances between the original data up to a small error---an important feature for many computational tasks. In this application, the number of hyperplanes $m$ determines the degree of compression that is achieved. Second, the map $t\mapsto\sign(At+\tau)$ occurs as a simple and popular model in signal processing. In this context, $Ax$ models the analog measurements of a signal $x$ and the map $y\mapsto \sign(y+\tau)$ a coarse quantizer that converts $Ax$ into a digital signal. Property \eqref{eqn:deltaTessProp} plays an essential role in proving error bounds for algorithms that aim to reconstruct a signal $x$ from its quantized measurements $\sign(Ax+\tau)$. In this second application, $m$ corresponds to the number of measurements of the signal. In both applications, it is of significant interest to understand which $m$ is necessary and sufficient to guarantee that $f$ is a $\delta$-embedding of a given set $T$. We refer to \cite{DM18,PlV14} and the references therein for a more extensive discussion of these applications.\par   
In this work we focus on what is arguably the most natural random hyperplane tessellation---with directions selected according to the standard gaussian distribution. Motivated by the aforementioned applications, our goal is to derive sharp bounds on the number of such hyperplanes that are needed to generate a $\delta$-tessellation. This problem was first studied by Plan and Vershynin \cite{PlV14} for subsets of the Euclidean sphere $S^{n-1}$, and with hyperplanes that pass through the origin (i.e., $\tau_i=0$ for $i=1,\ldots,m$).
\begin{Remark}
It is important to stress that hyperplanes that pass through the origin cannot separate points on any ray $\{ \alpha x : \alpha \geq 0\}$. As a result, shifts must be used in a $\delta$-tessellation of an arbitrary set.
\end{Remark}

To formulate the result from \cite{PlV14}, denote by $G$ the standard gaussian random vector and set
$$
\ell_*(T):= \E \sup_{t \in T} |\inr{G,t}| 
$$
denote the \emph{gaussian mean width} of a set $T\subset \R^n$. Also, let
$$
d_{S^{n-1}}(x,y)=\frac{1}{\pi}\arccos\left(\frac{\inr{x,y}}{\|x\|_2\|y\|_2}\right)
$$
be the spherical distance.

\begin{Theorem} \label{thm:PV} \cite{PlV14}
There are absolute constants $c_1$ and $c_2$ such that the following holds. Let $T \subset S^{n-1}$ and let $A:\R^n \to \R^m$ be the gaussian matrix. For $\delta>0$, if
$$
m \geq c_1 \frac{\ell_*^2(T)}{\delta^6}
$$
then with probability at least $1-2e^{-c_2m\delta^2}$, for every $x,y\in T$,
\begin{equation}
\label{eqn:PlVEmbed}
\left|\frac{1}{m}d_H(\sign(Ax),\sign(Ay)) - d_{S^{n-1}}(x,y)\right|\leq \delta.
\end{equation}
\end{Theorem}
\begin{Remark}
It was conjectured by Plan and Vershynin that the optimal bound on $m$ should be $m\sim \ell_*^2(T)/\delta^2$.
\end{Remark}
An improvement on the estimate from Theorem \ref{thm:PV} can be found in \cite{OyR15} (see Theorem~2.5 there). To formulate it, let ${\cal N} (T,\eps)$ be the covering number at scale $\eps$ with respect to the Euclidean norm; that is, the smallest number of open Euclidean balls of radius $\eps$ needed to cover $T$. Moreover, set $T-T=\{t_1-t_2 : t_1,t_2 \in T\}$.
\begin{Theorem} \cite{OyR15}
\label{thm:OyR}
There are absolute constants $c_1,c_2,$ and $c_3$ such that the following holds. Let $T \subset S^{n-1}$ and let $A:\R^n \to \R^m$ be the gaussian matrix. For $\delta>0$ let
$$
0<\theta \leq c_1 \frac{\delta}{\sqrt{\log(e/\delta)}}
$$
and set
\begin{equation} \label{eq:est-OyR15}
m \geq c_2 \left(\frac{\log{\cal N}(T,\theta)}{\delta^2} + \frac{\ell_*^2((T-T)\cap \theta B_2^n)}{\delta^3}\right).
\end{equation}
Then with probability at least $1-2e^{-c_3 m\delta^2}$, for every $x,y\in T$,
$$
\left|\frac{1}{m}d_H(\sign(Ax),\sign(Ay)) - d_{S^{n-1}}(x,y)\right|\leq \delta.
$$
\end{Theorem}

As noted previously, the hyperplanes in Theorems~\ref{thm:PV} and \ref{thm:OyR} cannot be used to generate a $\delta$-tessellation of an arbitrary set  $\R^n$. A natural way of addressing that issue is to consider parallel shifts---meaning that $\tau_1,\ldots,\tau_m$ need not be $0$. As it happens, under relatively mild randomness assumptions on $A$, \emph{random shifts} prove to be effective---leading to hyperplane tessellations that endow multiplicative \emph{isomorphic} embeddings in $\{-1,1\}^m$.

To formulate that fact, recall that a random vector $X$ is isotropic if for every $t \in \R^n$, $\E \inr{X,t}^2 = \|t\|_2^2$. The random vector is $L$-subgaussian if for every $p \geq 2$ and $t\in \R^n$,
$$
(\E |\langle X,t\rangle|^p )^{1/p} \leq L \sqrt{p} (\E \langle X,t\rangle^2)^{1/2}.
$$
Let $A$ be a matrix whose rows are independent copies of a symmetric, isotropic, $L$-subgaussian random vector and let $\tau_i$ be independent random variables, which are distributed uniformly in $[-\lambda,\lambda]$ and are independent of $A$.

\begin{Theorem} \label{thm:DM} \cite{DM18}
There are constants $c_0,...,c_4$ that depend only on $L$ such that the following holds. Let $T \subset \R^n$ and set $R=\sup_{t \in T} \|t\|_2$. Put $\lambda = c_0 R$ and set
\begin{equation} \label{eqn:bitCompDM18}
m \geq c_1R \log(eR/\delta) \frac{\ell_*^2(T)}{\delta^3}.
\end{equation}
Then with probability at least $1-8\exp(-c_2m\delta/R)$, for any $x,y \in {\rm conv}(T)$ that satisfy $\|x-y\|_2 \geq \delta$, one has
\begin{equation} \label{eq:isomorphic}
c_3\frac{ \|x-y\|_2}{R} \leq \frac{1}{m}d_H(\sign(Ax+\tau),\sign(Ay+\tau)) \leq c_4\sqrt{\log(eR/\delta)} \cdot \frac{\|x-y\|_2}{R}.
\end{equation}
\end{Theorem}

It is unrealistic to hope for an almost isometric estimate in the context of Theorem \ref{thm:DM} when using an arbitrary subgaussian matrix $A$ (because of the behaviour of the means $\E d_H(\sign(Ax+\tau),\sign(Ay+\tau))$. But as it happens, our first result is that if $A:\R^n \to \R^m$ is a gaussian matrix and $m$ is well-chosen, then the random mapping $t \mapsto \sign(At+\tau)$ is an almost isometric embedding (in the additive error sense) of an arbitrary $T \subset \R^n$. Moreover, the estimate on $m$ is of a similar nature to the one from \eqref{eq:est-OyR15}.

\begin{Theorem} \label{thm:main-A}
There exist absolute constants $c_0,...,c_3$ such that the following holds. Let $T \subset \R^n$ and put $R=\sup_{t \in T} \|t\|_2$. Set $\delta \in (0,\frac{R}{2}]$, $u\geq 1$ and let
$$
0<\theta \leq c_0 \frac{\delta}{\sqrt{\log(e\lambda/\delta)}}.
$$
Consider $\lambda \geq c_1 R\sqrt{\log(R/\delta)}$ and
\begin{equation}
\label{eqn:main-ABdm}
m\geq c_2 \left( \lambda^2 \frac{\log {\cal N}(T,\theta)}{\delta^2}  +  \lambda \frac{\ell_*^2((T-T)\cap \theta B^n_2)}{\delta^3} \right).
\end{equation}
If $A:\R^m \to \R^n$ is the standard gaussian matrix and $\tau$ is uniformly distributed in $[-\lambda,\lambda]^m$, then with probability at least $1-2\exp(-c_3\delta^2m/\lambda^2)$, the map $f(t)= \sign(At+\tau)$ satisfies
\begin{equation} \label{eqn:GaussianDitheredText}
\sup_{x,y\in T}\left|\frac{\sqrt{2\pi}\lambda}{m}d_H(f(x),f(y))-\|x-y\|_2\right|\leq  \delta.
\end{equation}
\end{Theorem}

The proof of Theorem \ref{thm:main-A} follows from a generic embedding result, presented in Section~\ref{sec:generic}. We show that if a deterministic matrix $A:\R^n \to \R^m$ `acts well' on $T$ (in a sense that is clarified in what follows) and $\tau$ is as in Theorem \ref{thm:main-A}, then with high probability with respect to $\tau$, the mapping $t \mapsto \sign(At+\tau)$ is a $\delta$-embedding of $T$. We then show in Section \ref{sec:gaussian} that for any $T \subset \R^n$, a typical realization of the gaussian matrix $A$ `acts well' on $T$.

\vskip0.3cm

The choice of $\lambda$ is very natural --- random shifts at a scale proportional to the radius of the set ensure that one can separate points that belong to the same ray $\{\alpha x : \alpha \geq 0\}$. In fact, $x \mapsto \sign(Ax+\tau)$ can already fail to be a $\delta$-embedding with high probability on sets consisting of only two points on a ray if $\lambda$ is smaller than the radius (see Appendix~\ref{sec:minShift} for details). At the same time, the estimate \eqref{eqn:main-ABdm} in Theorem \ref{thm:main-A} looks anything but natural. In fact, it is reasonable to guess that it is suboptimal. It was conjectured that the right tradeoff between the accuracy parameter $\delta$, the set $T$ and the dimension $m$ is $m \sim \ell_*^2(T)/\delta^2$---just as in the conjecture regarding subsets of the sphere.
\begin{tcolorbox}
Surprisingly, the assertion of Theorem \ref{thm:main-A} is optimal (up to logarithmic factors). And with the natural choice of $\lambda$, the conjectured optimal estimate of $m \sim \ell_*^2(T)/\delta^2$ is simply false.
\end{tcolorbox}
We prove that fact by establishing two lower bounds, the first of which deals with properties of arbitrary random embeddings into the discrete cube.

\begin{Definition} \label{def:random-emb}
Let $(f,\phi)$ be a pair of maps, where $f:B_2^n \to \{-1,1\}^m$ and $\phi:\{-1,1\}^m \to \R$. $(f,\phi)$ is a \emph{$\delta$-inner product-preserving embedding} of $W \subset B_2^n$ if 
\begin{equation}
\label{eqn:defInnProdPres}
\sup_{(x,y) \in W \times W} \left| \phi(f(x),f(y)) - \langle x,y\rangle \right| < \delta.
\end{equation}
A pair of random maps $(f,\phi)$ is a \emph{random inner product-preserving embedding into $\{-1,1\}^m$ with parameters $n,\delta,\eta,$ and $N$} if, for any $W \subset B_2^n$ of cardinality $N$, \eqref{eqn:defInnProdPres} holds with probability at least $1-\eta$.
\end{Definition}
In what follows we abuse notation by omitting the dependence of the embedding on the parameters $n$, $\eta$ and $N$ and keeping track only of the accuracy parameter $\delta$ and the embedding dimension $m$. Using a polarization argument, it is straightforward to verify that a $\delta$-embedding of $W\cup(-W)$ yields a $\delta$-inner product-preserving embedding of $W$. In particular, the map in Theorem~\ref{thm:main-A} yields a random inner product-preserving embedding in the sense of Definition \ref{def:random-emb} if
\begin{equation} \label{eq:est-m-rand-emb}
m\sim \log(e/\delta)\frac{\log(eN/\eta)}{\delta^2}.
\end{equation}
The first lower bound we present shows that \eqref{eq:est-m-rand-emb} is not improvable (up to a logarithmic factor).
\begin{Theorem} \label{thm:main-B-1}
There are absolute constants $c_0,c_1,c_2,c_3$ such that the following holds. Let $N$ be an integer and set $0<\eta<1$. Let $c_0(\eta/N)^{1/2} \leq \delta \leq c_1$ and put $n=\frac{c_2}{\delta^2}\log(eN/\eta)$. If $(f,\phi)$ is a random inner product-preserving embedding into $\{-1,1\}^m$ with parameters $n,\delta,\eta,N$ then
$$
m \geq c_3 \frac{\log (eN/\eta)}{\delta^2}.
$$
\end{Theorem}

The proof of Theorem \ref{thm:main-B-1} actually reveals more than is stated: for any random embedding as above, there is a `bad set' $W$ that forces the embedding dimension to be at least $\delta^{-2} \log (N/\eta)$. Moreover, there is an absolute constant $\alpha$ so that $\hat{W}=W\cup(-W)$ is $\alpha$-separated. In particular, if $\delta$ is sufficiently small, then $(\hat{W}-\hat{W}) \cap \delta B_2^n = \{0\}$. As a result, for this bad set the dominating term in the upper estimate on the embedding dimension in Theorem \ref{thm:main-A} is the entropic one, showing that it cannot be improved (up to logarithmic factors)---regardless of the choice of the random embedding.

The proof, which can be found in Section \ref{sec:proof-B1}, uses similar ideas to the ones appearing in a construction due to Alon and Klartag \cite{AlK17}. Alon and Klartag studied $\delta$-distance sketches, which are data structures that allow the reconstruction of all Euclidean scalar products between points in a subset of $\R^n$ of a given cardinality.

\vskip0.3cm

In the second lower bound, our focus is on the specific embedding $t \mapsto \sign(At+\tau)$ that is used in Theorem \ref{thm:main-A}. We show that for any convex body\footnote{A convex body in $\R^n$ is a convex, centrally-symmetric set with a nonempty interior.} $T \subset \R^n$, the best possible dimension $m$ one can hope for coincides with the estimate from Theorem \ref{thm:main-A} (up to logarithmic factors). To formulate that lower bound, set
$$
d^*(T) = \left(\frac{\ell_*(T)}{\sup_{t \in T}\|t\|_2}\right)^2
$$
to be the \emph{Dvoretzky-Milman dimension} of $T$.
\begin{Theorem} \label{thm:main-B-2}
There are absolute constants $c_0$, $c_1$ and $c_2$ such that the following holds. For any convex body $T \subset \R^n$, if $\lambda \geq c_0 \ell_*((T-T)\cap \delta B^n_2)/\sqrt{m}$,
$$
m \leq c_1 \lambda \frac{\ell_*^2((T-T)\cap \delta B^n_2)}{\delta^3},
$$
and $f$ is as in Theorem~\ref{thm:main-A}, then with probability at least $1-2\exp(-c_2d^*(T))$ there are $x,y \in T$ such that
$$
\left|\frac{\sqrt{2\pi}\lambda}{m}d_H(f(x),f(y))-\|x-y\|_2\right| \geq 2\delta.
$$
\end{Theorem}

\begin{Remark}
The lower bound on $\lambda$ in Theorem~\ref{thm:main-B-2} is not really restrictive. As an example, we will present in Lemma~\ref{lemma:lower-on-lambda} a rather general situation in which that lower bound holds. 
\end{Remark}

There is a purely geometric reason for Theorem \ref{thm:main-B-2} being true. We first show that it is hard to embed a Euclidean ball using the mapping $y \mapsto \sign(y + \tau)$. More accurately, if there is a subset of coordinates $I \subset \{1,...,m\}$ such that $P_I V=\{ (x_i)_{i \in I} : x \in V\} =r B_2^I$, then embedding $V$ in $\{-1,1\}^m$ by using the mapping $y \mapsto \sign(y+\tau)$ fails unless $m$ is sufficiently large (depending on $|I|$, $r$ and $\lambda$). We then invoke the classical Dvoretzky-Milman Theorem and show that if $T \subset \R^n$ is a convex body and $A:\R^n \to \R^m$ is the standard gaussian matrix, then for a well-chosen $I$, with high probability, $P_I V=P_I A (T \cap \delta B_2^n)$ contains a large Euclidean ball. In particular, with high probability, the map $t \mapsto \sign(At + \tau)$ fails to be a $\delta$-embedding of $T$ unless $m$ is as in Theorem \ref{thm:main-B-2}.

The proof of Theorem \ref{thm:main-B-2} can be found in Section \ref{sec:proof-B2}. The rest of that section is devoted to the conjecture we mentioned previously: that $m \sim \ell_*^2(T)/\delta^2$ suffices to ensure that $t \mapsto \sign(At + \tau)$ is a $\delta$-embedding of an arbitrary $T \subset \R^n$. Thanks to Theorem~\ref{thm:main-B-2}, we construct an example which shows that the conjecture is false: a set for which $m \sim \ell_*^2(T)/\delta^3$ is required---see Theorem~\ref{thm:counter-example}. 

\section{A generic embedding result} \label{sec:generic}
The first step in the proof of Theorem \ref{thm:main-A} is actually rather general. We consider a set $T \subset \R^n$, and a (deterministic) matrix $A:\R^n \to \R^m$. For a parameter $\lambda>0$, let $\tau$ be a random vector distributed uniformly  in $[-\lambda,\lambda]^m$, and the embedding function we consider is $f: \R^n \to \{-1,1\}^m$, defined by
\begin{equation} \label{eqn:GaussBinEmdDef}
 f(x)=\sign (Ax+\tau).
\end{equation}
Here, as always, the sign-function is applied component-wise. 

Our goal is to show that if $A$ acts on $T$ `in a regular way', then $f$ is a $\delta$-embedding of $T$ into $\{-1,1\}^m$. To explain what we mean by `regularity', let $0 < \theta < \delta$, and set $T_\theta \subset T$ to be a minimal $\theta$-net of $T$; in particular, $|T_\theta| = {\cal N}(T,\theta)$. For $v \in \R^m$ and $1 \leq s \leq m$ set
$$
\|v\|_{[s]} = \max_{|I|=s} \left(\sum_{i \in I} v_i^2 \right)^{1/2},
$$
and assume that the following holds:
\begin{description}
\item{$(a)$} \underline{uniform $\ell_1$-concentration:} there is a constant $\kappa$ such that
\begin{equation} \label{eqn:ell1ell2assump}
\sup_{x,y\in T_{\theta}}\left|\frac{\kappa}{m} \|A(x-y)\|_1 - \|x-y\|_2\right|\leq \delta.
\end{equation}
\item{$(b)$} \underline{$A$ maps $T$ to `well-spread' vectors:} For $s=\lfloor\delta m/\lambda\rfloor$,
\begin{equation} \label{eqn:knormbias}
\frac{1}{\sqrt{s}} \sup_{x\in T_\theta} \|Ax\|_{[s]}\leq \lambda,
\end{equation}
and
\begin{equation} \label{eqn:knormoscillations}
\frac{1}{\sqrt{s}} \sup_{x\in (T-T)\cap \theta B_2^n} \|Ax\|_{[s]}\leq \delta.
\end{equation}
\end{description}

Denote the normalized Hamming distance on $\{-1,1\}^m$ by
\begin{equation*}
\tilde{d}(x,y)= \frac{2\lambda \kappa}{m}  d_H(x,y)
\end{equation*}
for a constant $\kappa$ whose value is specified in what follows. In our application, where $A$ is the gaussian matrix, $\kappa$ turns out to be an absolute constant. 
\vskip0.3cm
Once $A$ exhibits the necessary regular behaviour, one may show that with high probability with respect to $\tau$, $f$ is a $\delta$-embedding of $T$.

\begin{Theorem} \label{thm:mainGeneric}
There exist absolute constants $c_1,c_2,$ and $c_3$ such that the following holds. Let $\theta \leq \delta$, set
$$
m \geq c_1 \lambda^2 \kappa^2\frac{\log{\cal N}(T,\theta)}{\delta^2}
$$
and assume that $A$ satisfies \eqref{eqn:ell1ell2assump}, \eqref{eqn:knormbias}, and \eqref{eqn:knormoscillations}. Then with probability at least
$$
1-2 \exp(-c_2 \delta^2 m/(\lambda^2 \kappa^2))
$$
with respect to $\tau$,
\begin{equation*}
\sup_{x,y\in T}\left|\tilde{d}(f(x),f(y))-\|x-y\|_2\right|\leq c_3 (\kappa+1) \delta.
\end{equation*}
\end{Theorem}

The proof of Theorem \ref{thm:mainGeneric} requires two simple observations.
\begin{Lemma}\label{lem:exp}
Fix $\lambda>0$ and let $\tau$ be uniformly distributed in $[-\lambda, \lambda]$. If $\phi_\lambda(x)=(|x|-\lambda)\IND_{\{|x|\geq \lambda\}}$ then for every $a,b\in \R$,
\begin{equation*}
\left|2\lambda \PP(\sign(a+\tau) \neq \sign(b+\tau))- |a-b|\right|\leq \phi_\lambda(a)+\phi_\lambda(b).
\end{equation*}
\end{Lemma}

The proof of Lemma \ref{lem:exp} follows from a straightforward and tedious computation. It is deferred to Appendix \ref{app:proof-lem:exp}.

\begin{Corollary}  \label{cor:exp}
Let $x,y \in T_\theta$. If \eqref{eqn:knormbias} holds then
\begin{equation*}
\left| \E_\tau d_H (f(x),f(y)) - \frac{1}{2\lambda} \|A(x-y)\|_1 \right| \leq s.
\end{equation*}
\end{Corollary}

\proof Observe that
\begin{align*}
\E_\tau d_H (f(x),f(y)) & =\E_\tau \sum_{i=1}^m \IND_{\{ \sign( (Ax)_i + \tau_i) \not= \sign( (Ay)_i + \tau_i)\}}
\\
& = \sum_{i=1}^m \PP_\tau \left(\sign( (Ax)_i + \tau_i) \not= \sign( (Ay)_i + \tau_i)\right),
\end{align*}
and therefore, by Lemma \ref{lem:exp},
$$
\left| \E_\tau d_H (f(x),f(y)) - \frac{1}{2\lambda} \|A(x-y)\|_1 \right|  \leq \frac{1}{\lambda} \sup_{z \in T_\theta}   \sum_{i=1}^m |(Az)_i| \IND_{\{|(Az)_i|>\lambda\}}.
$$
To control the right-hand side, note that by \eqref{eqn:knormbias}, for every $z \in T_\theta$,
\begin{equation*} 
|(Az)^*_{s}| \leq \frac{1}{\sqrt{s}} \|Az\|_{[s]}\leq \lambda,
\end{equation*}
where $(Az)^*_{s}$ denotes the $s$-largest coordinate of $( |(Az)_i|)_{i=1}^m$. In particular,
\begin{equation} \label{eqn:largeCoeffk-normBias}
\sup_{z\in T_{\theta}} \left| \left\{i \in \{1,...,m\} \ : \ |(Az)_i|> \lambda \right\} \right| \leq s.
\end{equation}
Thus, using \eqref{eqn:knormbias} once again,
$$
\frac{1}{\lambda}\sup_{z \in T_\theta} \sum_{i=1}^m |(Az)_i| \IND_{\{|(Az)_i|>\lambda\}} \leq \frac{1}{\lambda}\sup_{z \in T_\theta} \sqrt{s} \|Az\|_{[s]} \leq s.
$$
\endproof

\noindent{\bf Proof of Theorem \ref{thm:mainGeneric}.} For every $x\in T$ let $\pi x  \in T_\theta$ satisfy that $\|x-\pi x\|_2 \leq \theta$. By the triangle inequality,
\begin{align}\label{eq:four_summands}
& \left|\tilde{d}(f(x),f(y))-\|x-y\|_2\right| \nonumber\\
& \qquad \leq   \left|\tilde{d}(f(x),f(y))- \tilde{d}(f(\pi x),f(\pi y))\right| \nonumber
\\
& \qquad \qquad + \left|\tilde{d}(f(\pi x),f(\pi y)) - \E_{\tau} \tilde{d}(f(\pi x),f(\pi y)) \right|
+ \left|\E_{\tau} \tilde{d}(f(\pi x),f(\pi y)) - \frac{\kappa}{m}\|A(\pi x - \pi y)\|_1\right| \nonumber
\\
& \qquad \qquad + \left|\frac{\kappa}{m}\|A(\pi x - \pi y)\|_1 - \|\pi x-\pi y\|_2\right|
+ \left| \|\pi x-\pi y\|_2-\|x-y\|_2\right| \nonumber
\\
& \qquad = (a)+(b)+(c)+(d)+(e).
\end{align}
Clearly, $(e)\leq 2\theta \leq 2\delta$ and, by \eqref{eqn:ell1ell2assump}, $(d)\leq \delta$. Moreover, it follows from Corollary \ref{cor:exp} that
\begin{align*}
(c)= & \left|\E_{\tau} \tilde{d}(f(\pi x),f(\pi y)) - \frac{\kappa}{m}\|A(\pi x - \pi y)\|_1\right|
\\
= & \frac{2\lambda \kappa}{m} \left|\E_{\tau} d_H(f(\pi x),f(\pi y)) - \frac{1}{2\lambda}\|A(\pi x - \pi y)\|_1\right| \leq 2\lambda \kappa \frac{s}{m} \leq 2 \kappa \delta.
\end{align*}

Turning to $(b)$, note that
$$
\tilde{d}(f(x),f(y))=\frac{2\lambda \kappa}{m} \sum_{i=1}^m \IND_{\{\sign(Ax+\tau)_i \not = \sign(Ay+\tau)_i\}},
$$
and in particular it is the sum of independent, $\{0,1\}$-valued random variables. Hence, by Hoeffding's inequality there is an absolute constant $c$ such that,  for every $x,y\in T_\theta$ and $\delta>0$,
\begin{equation*}
\PP_{\tau}\left(\left|\tilde{d}(f(\pi x),f(\pi y)) - \E_{\tau} \tilde{d}(f(\pi x),f(\pi y))\right| \geq \delta \right)\leq 2 \exp\left(-c\frac{\delta^2 m}{\lambda^2 \kappa^2}\right).
\end{equation*}
It follows that if
$$
m \geq c_1 \kappa^2 \frac{\lambda^2}{\delta^{2}} \log(|T_{\theta}|),
$$
then by the union bound, with probability at least $1-2 \exp(-c_2\delta^2 m/(\lambda^2\kappa^2))$,
\begin{equation*}
\sup_{x,y\in T} \left|\tilde{d}(f(\pi x),f(\pi y)) - \E_{\tau} \tilde{d}(f(\pi x),f(\pi y))\right|\leq \delta.
\end{equation*}
Finally, to control $(a)$, by the triangle inequality we have that
\begin{align*}
& \left|\tilde{d}(f(x),f(y))- \tilde{d}(f(\pi x),f(\pi y))\right|
\\
& \qquad \leq \frac{2\lambda \kappa}{m} \sum_{i=1}^m \left|\IND_{\{ \sign((Ax)_i+\tau_i) \neq \sign ((Ay)_i+\tau_i)\}} -\IND_{ \{\sign((A \pi x)_i+\tau_i)\neq \sign((A \pi y)_i+\tau_i)\}}\right|.
\end{align*}
Observe that if
$$
\sign((Ax)_i+\tau_i)=\sign((A \pi x)_i+\tau_i) \ \ {\rm and} \ \
\sign((Ay)_i+\tau_i)=\sign((A\pi y)_i+\tau_i),
$$
then
$$
\IND_{\{\sign( (Ax)_i+\tau_i) \neq \sign((Ay)_i+\tau_i)\}} -\IND_{\{ \sign ((A\pi x)_i+\tau_i) \neq \sign( (A \pi y)_i+\tau_i)\}}=0.
$$
Therefore,
\begin{equation}
\label{eq:applyDM18}
\sup_{x,y\in T}\left|\tilde{d}(f(x),f(y))- \tilde{d}(f(\pi x),f(\pi y))\right|
\leq \frac{4\lambda \kappa}{m} \sup_{x \in T} \sum_{i=1}^m \IND_{\{\sign( (Ax)_i+\tau_i)\neq \sign( (A \pi x)_i+\tau_i)\}}.
\end{equation}
Moreover,
$$
\IND_{\{\sign( (Ax)_i+\tau_i)\neq \sign( (A \pi x)_i+\tau_i)\}} = \IND_{\{\tau_i \in [-(Ax)_i,-(A \pi x)_i] \cup [-(A \pi x)_i,-(Ax)_i]\}}
$$
implying that
\begin{equation} \label{eqn:indZeroSep}
\IND_{\{\sign((Ax)_i+\tau_i)\neq \sign( (A \pi x)_i+\tau_i)\}}=0
\end{equation}
on the set
$$
A_{\delta,i} = \left\{x\in \R^n \ : \ |(A \pi x)_i+\tau_i|>\delta \geq |(A (x-\pi x))_i|\right\}.
$$
Therefore,
\begin{equation} \label{eqn:signChangeCov}
\IND_{\{\sign((Ax)_i+\tau_i) \neq \sign( (A \pi x)_i+\tau_i)} \leq \IND_{A_{\delta,i}^c} \leq \IND_{\{|(A \pi x)_i+\tau_i|\leq \delta\}} + \IND_{\{| (A (x-\pi x))_i|>\delta\}},
\end{equation}
and in particular,
\begin{align}\label{eq:proof:Gaussian:twoprocesses}
&\sup_{x,y\in T}\left|\tilde{d}(f(x),f(y))- \tilde{d}(f(\pi x),f(\pi y))\right| \nonumber
\\
& \qquad \leq \frac{4\lambda \kappa}{m} \left(\sup_{y\in T_{\theta}} \sum_{i=1}^m \IND_{\{|(Ay)_i+\tau_i|\leq \delta\}}
  + \sup_{z\in (T-T)\cap \theta B_2^n} \sum_{i=1}^m \IND_{\{|(Az)_i|>\delta\}}\right).
\end{align}
Clearly, for every $\alpha \in \R$ and $1 \leq i \leq m$, $\PP_\tau(|\alpha+\tau_i| \leq \delta) \leq \delta/\lambda$; thus,
$$
\PP_{\tau}\left(| (Ay)_i+\tau_i|\leq \delta\right)\leq \frac{\delta}{\lambda}.
$$
By the Chernoff bound, there is an absolute constant $c_3$ such that with $\tau$-probability at least $1-\exp(-c_3m\delta/\lambda)$,
$$
\sum_{i=1}^m \IND_{\{|(Ay)_i+\tau_i|\leq \delta\}}\leq \frac{2\delta m}{\lambda}.
$$
Hence, if
$$
m \geq c_4 \frac{\lambda}{\delta} \log|T_\theta|
$$
then by the union bound, with $\tau$-probability at least $1-\exp(-c_3m\delta/(2\lambda))$
$$
\sup_{y\in T_{\theta}} \frac{4\lambda \kappa}{m}\sum_{i=1}^m \IND_{\{|(Ay)_i+\tau_i|\leq \delta\}} \leq 8\kappa \delta.
$$
Finally, by \eqref{eqn:knormoscillations}, for every $z\in (T-T)\cap \theta B_2^n$,
\begin{equation*}
(Az)^*_{s} \leq \frac{1}{\sqrt{s}}\|Az\|_{[s]}\leq \delta,
\end{equation*}
implying that the second term on the right-hand side of \eqref{eq:proof:Gaussian:twoprocesses} is bounded by $c\kappa \delta$.
\endproof

\subsection{The gaussian case} \label{sec:gaussian}

When $A$ is the standard gaussian matrix, the application of Theorem~\ref{thm:mainGeneric} is straightforward. The analysis is based on two well-known facts on the way that a gaussian matrix acts on an arbitrary subset of $\R^n$. The first is a consequence of the standard gaussian concentration inequality, see e.g.\ \cite[Lemma 2.1]{PlV14}.
\begin{Theorem} \label{thm:gaussian-1}
There is an absolute constant $c$ such that the following holds. Let $T \subset \R^n$ and $A:\R^n \to \R^m$ be the standard gaussian matrix. Then for $u \geq 1$,
$$
\PP \left( \sup_{t \in T} \left| \|A t\|_1 - \E \|At\|_1 \right| \geq u\sqrt{m} \ell_*(T) \right) \leq 2\exp\left(-cu^2 d^*(T)\right),
$$
where
$$
d^*(T) = \left(\frac{\ell_*(T)}{{\cal R}(T)}\right)^2 \ \ {\rm and} \ \ {\cal R}(T)=\sup_{t \in T} \|t\|_2.
$$
\end{Theorem}
The second fact required here is that a typical realization of $A$ maps $T$ to a set of `regular vectors'.
\begin{Theorem} \label{thm:gaussian-good-position}
There are absolute constants $c$ and $C$ such that the following holds. Let $T \subset \R^n$ and $A:\R^n \to \R^m$ be the standard gaussian matrix. Then for $u \geq 1$ and $1 \leq k \leq m$, with probability at least $1-2\exp(-cu^2 k\log(em/k))$,
$$
\sup_{t \in T} \|At\|_{[k]} \leq C \left(\ell_*(T) + u {\cal R}(T) \sqrt{k \log (em/k)} \right).
$$
\end{Theorem}
Theorem~\ref{thm:gaussian-good-position} follows by a standard chaining argument. We include this proof in Appendix~\ref{sec:knormGaussian} for the sake of completeness. Using Theorem \ref{thm:gaussian-1} and Theorem \ref{thm:gaussian-good-position}, it is straightforward to establish Theorem \ref{thm:main-A}.

\noindent{\bf Proof of Theorem \ref{thm:main-A}.} Following Theorem \ref{thm:mainGeneric}, we apply Theorem \ref{thm:gaussian-good-position} for the sets $T_{\theta}$ and $(T-T) \cap \theta B_2^n$. Thus, for every fixed $1 \leq k \leq m$, with probability at least $1-2\exp(-c_0 u^2 k \log(em/k))$, for every $t \in T_{\theta}$,
\begin{equation} \label{eq:gaussian-k-norm-1}
\|At\|_{[k]} \leq c_1\left(\ell_*(T_{\theta})+ u {\cal R}(T_{\theta}) \sqrt{k \log (em/k)} \right)
\end{equation}
and for every $z \in (T-T) \cap \theta B_2^n$,
\begin{equation} \label{eq:gaussian-k-norm-2}
\|At\|_{[k]} \leq c_1\left(\ell_*((T-T)\cap \theta B_2^n)+ u \theta \sqrt{k \log (em/k)} \right).
\end{equation}
Using these inequalities for $k=\lfloor\delta m/\lambda\rfloor$, it is straightforward to verify that \eqref{eqn:knormbias} and \eqref{eqn:knormoscillations} hold under the stated conditions and with the stated probability. Moreover, Theorem~\ref{thm:gaussian-1} immediately implies that \eqref{eqn:ell1ell2assump} holds with the wanted probability for $\kappa=\sqrt{\frac{\pi}{2}}$. The result now follows from Theorem~\ref{thm:mainGeneric}.

\endproof

\begin{tcolorbox}
Interesting as Theorem \ref{thm:main-A} may be, it is not really a major surprise that the embedding $t \mapsto \sign(At+\tau)$ is successful. What is far more surprising is that the seemingly unnatural estimate on the dimension happens to be sharp. The next two sections are devoted to the proof of that fact.
\end{tcolorbox}

\section{Proof of Theorem \ref{thm:main-B-1}} \label{sec:proof-B1}

The proof is based on two standard observations (see, e.g., \cite{BLM13} for their proofs).
\begin{Lemma} \label{lemma:standard}
Let $X$ be uniformly distributed in $S^{n-1}$. Then for $0 < \delta <1/2$ and any $z \in S^{n-1}$,
$$
\PP( |\langle X,z\rangle| \geq \delta  ) \geq \exp(-2\delta^2 n),
$$
and
$$
\PP( \|X-z\|_2 \leq \alpha) \leq \frac{1}{\sqrt{n}} \alpha^{n-1}.
$$
\end{Lemma}

Set $T \subset S^{n-1}$ to be a $1/2$-separated subset of $S^{n-1}$ of cardinality  $|T| = \exp(c_0n)$. Such a set exists thanks to a standard volumetric argument. The value of the constant $c_0$ remains unchanged for the rest of this section.

\vskip0.3cm

\begin{Lemma} \label{lemma:first-obs}
There is an absolute constant $c_1$ such that the following holds. Let
\begin{equation} \label{eq:cond-on-k-1}
\log k \geq c_1 \max\{\delta^2 n, \log n\}
\end{equation}
and let $X_1,...,X_k$ be independent and uniformly distributed in $S^{n-1}$. Then with probability at least $0.99$, for every $t_1,t_2 \in T$, there is a $1 \leq j \leq k$ such that
$$
|\langle X_j,t_1-t_2\rangle| \geq \delta.
$$
\end{Lemma}
\proof Define $U\subset S^{n-1}$ by
$$
U=\left\{\frac{t_i-t_j}{\|t_i-t_j\|_2} : t_i \not = t_j, \ t_i,t_j \in T \right\},
$$
and let $\gamma=2\delta$. 
Since $\|t_i-t_j\|_2 \geq 1/2$ if $t_i \not = t_j$, it suffices to show that, with the stated probability, for any $u\in U$ there is some $1 \leq j \leq k$ such that
$$
|\langle X_j,u\rangle| \geq \gamma.
$$
Let $X$ be uniformly distributed in $S^{n-1}$ and set $p=\bP(|\langle X,u\rangle| \geq \gamma)$. By independence of the $X_j$,
$$
\PP( \forall 1 \leq j \leq k, : \ |\langle X_j,u\rangle| < \gamma) = [\bP(|\langle X,u\rangle| < \gamma)]^k = (1-p)^k\leq e^{-kp}.
$$
By the union bound, the wanted estimate follows if $|T|e^{-kp}\leq 0.01$. Recalling that $|T| =\exp(c_0n)$ and that, by Lemma~\ref{lemma:standard}, $p \geq \exp(-c_2\gamma^2 n)$, the condition holds when $k$ satisfies \eqref{eq:cond-on-k-1}.
\endproof

The second observation we require follows immediately from the second part of Lemma~\ref{lemma:standard} and the union bound.
\begin{Lemma} \label{lemma:second-obs}
There are absolute constants $c_1$ and $c_2$ such that the following holds. Let $0 < \alpha \leq c_1$ and $\log k \leq c_2 n \log(1/\alpha)$. Then with probability at least $0.99$, for every $1 \leq i \not = j \leq k$,
$\|X_i \pm X_j\|_2 \geq \alpha$ and for every $1 \leq i \leq k$ and every $t \in T$, $\|X_i \pm t\|_2 \geq \alpha$.
\end{Lemma}

Combining Lemma \ref{lemma:first-obs} and Lemma \ref{lemma:second-obs}, we find that there are (deterministic) sets $V,T \subset S^{n-1}$ with the following properties:
\begin{description}
\item{$(1)$} $V\cup(-V) \cup T$ is an $\alpha$-separated subset of $S^{n-1}$.
\item{$(2)$} For every $t_1,t_2 \in T$ there is some $v \in V$ such that $|\langle v,t_1-t_2\rangle| \geq \delta$.
\end{description}
Next, set $s =\frac{k}{N-1}$; split the set $V$ to $s$ disjoint sets, each  of cardinality $N-1$; and denote those sets by $V_j$, $1 \leq j \leq s$. Let $Y$ be uniformly distributed in $T$ and consider the random sets
$$
V_j \cup \{Y\}, \ \ 1 \leq j \leq s,
$$
which are of cardinality $N$. Note that the only source of randomness in $V_j \cup \{Y\}$ is $Y$. Let $(f,\phi)$ be a random inner product-preserving embedding into $\{-1,1\}^m$ with parameters $n,\frac{\delta}{4},\eta,$ and $N$ (see Definition~\ref{def:random-emb}), and consider the event
$$
\left\{ (f,\phi) \ \ \delta/4{\rm -embeds \ } V_j \cup \{Y\} \ {\rm in \ } \{-1,1\}^m \ {\rm for \ every \ } 1 \leq j \leq s \right\}.
$$
Observe that for every fixed realization of $Y$ and the union bound,
$$
\PP_{(f,\phi)} \left( \left\{ (f,\phi) \ \ \delta/4{\rm-embeds \ } V_j \cup \{Y\} \ {\rm in \ } \{-1,1\}^m \ {\rm for \ every \ } 1 \leq j \leq s \right\} \right) \geq 1-s\eta \geq 1/2
$$
provided that
\begin{equation} \label{eq:cond-on-k-11}
s = \frac{k}{N-1} \leq \frac{1}{2\eta}.
\end{equation}
In that case,
\begin{align*}
& \PP \left( \left\{ (f,\phi) \ \delta/4{\rm -embeds \ } V_j \cup \{Y\} \ {\rm in \ } \{-1,1\}^m \ {\rm for \ every \ } 1 \leq j \leq s \right\} \right)
\\
& \qquad = \E_Y \left( \PP_{(f,\phi)} \left( \left\{ (f,\phi) \ \delta/4{\rm -embeds \ } V_j \cup \{Y\} \ {\rm in \ } \{-1,1\}^m \ {\rm for \ every \ } 1 \leq j \leq s \right\} \right) \big| Y \right) \geq \frac{1}{2}.
\end{align*}
Hence, by Fubini's Theorem, there is some realization of $(f,\phi)$ such that
$$
\PP_Y \left( \left\{ (f,\phi) \ \delta/4{\rm -embeds \ } V_j \cup \{Y\} \ {\rm in \ } \{-1,1\}^m \ {\rm for \ every \ } 1 \leq j \leq s \right\} \right) \geq \frac{1}{2}.
$$
Thus, by the definition of $Y$, there is a set $T^\prime \subset T$ of cardinality at least $\frac{1}{2}|T|=\frac{1}{2}\exp(c_0n)$, such that for every $t \in T^\prime$ and every $1 \leq j \leq s$, the set $V_j \cup \{t\}$ is $\delta/4$-embedded in $\{-1,1\}^m$ by the fixed embedding $(f,\phi)$.

\vskip0.3cm
The key observation is now as follows.

\begin{Lemma} \label{lemma:lower-injection}
The embedding $f$ is injective on $T^\prime$.
\end{Lemma}

\proof  Let $t_1, t_2  \in T^\prime$ such that $f(t_1)=f(t_2)$.
By Property $(2)$, there is a $v \in V$ that `separates' $t_1$ and $t_2$; i.e.,
$$
\left|\langle t_1,v\rangle - \langle t_2,v\rangle\right| \geq \delta.
$$
Let $j$ be such that $v \in V_j$, observe that $(f,\phi)$ is a $\delta/4$-embedding of $V_j \cup \{t_1\}$ and $V_j \cup \{t_2\}$, and trivially we have that $\phi(f(t_1),f(v))=\phi(f(t_2),f(v))$. Because $(f,\phi)$ is a $\delta/4$-embedding of both sets, it follows that
\begin{equation} \label{eq:good-embedding}
\left| \phi(f(t_1),f(v)) - \langle t_1,v\rangle \right| \leq \frac{\delta}{4} \ \ {\rm and} \ \ \left| \phi(f(t_2),f(v)) - \langle t_2,v\rangle \right| \leq \frac{\delta}{4}.
\end{equation}
Hence,
$$
|\langle t_1,v\rangle - \langle t_2,v\rangle| \leq \left|\langle t_1,v\rangle-\phi(f(t_1),f(v)) \right| + \left|\langle t_2,v\rangle-\phi(f(t_2),f(v)) \right| \leq \frac{\delta}{2},
$$
which is a contradiction.
\endproof

Thanks to Lemma \ref{lemma:lower-injection} we have that $2^m \geq |T^\prime| = \frac{1}{2}\exp(c_0n)$; thus, $m \geq c_1 n$, and all that is left is to verify that it is possible to take $n$ sufficiently large. By \eqref{eq:cond-on-k-11}, one may set $k=(N-1)/2\eta \sim N/\eta$ and by Lemma~\ref{lemma:first-obs} one has to ensure that
$$
\log(c_1N/\eta) \geq c_2 \max\{\delta^2 n, \log n\}.
$$
Setting $n=c_3 \frac{\log(c_1N/\eta)}{\delta^2}$, the wanted condition holds by the choice of $\delta \geq c_4(\eta/N)^{1/2}$, and the condition in Lemma \ref{lemma:second-obs} is satisfied for a well-chosen absolute constant $\alpha$ provided that $\delta\leq c_5$.
\endproof

\section{Proof of Theorem \ref{thm:main-B-2}} \label{sec:proof-B2}

The proof is based on two ingredients. First, the Dvoretzky-Milman Theorem (for projections). For more information on the Dvoretzky-Milman Theorem and its central role in Asymptotic Geometry Analysis, see e.g.\ \cite{MR3331351}.
\begin{Theorem} \label{thm:DvorMil}
There are absolute constants $c_1,c_2,$ and $c_3$ such that the following holds. Let $K \subset \R^n$ be a convex body, and consider $s\leq c_1d^*(K)$. If $A:\R^n \to \R^s$ is the gaussian matrix, then with probability at least $1-2\exp(-c_2 d^*(K))$,
$$
c_3 \ell_*(K) B_2^s \subset A K.
$$
\end{Theorem}
The proof of Theorem \ref{thm:DvorMil} is standard. It can be found, for example, in \cite{MR3571258}.

\begin{Remark}
The Dvoretzky-Milman Theorem actually implies a two-sided, almost isometric equivalence (i.e., $c_3$ can be made arbitrarily close to $1$, and the lower bound is complemented by an analogous upper bound). However, for our purposes, an isomorphic, one-sided bound suffices.
\end{Remark}

The second component needed here is a standard probabilistic estimate.
\begin{Lemma} \label{lemma:lower-lambda}
There are absolute constants $c_0,c_1$ such that the following holds. Let $\tau=(\tau_i)_{i=1}^m$ be distributed uniformly in $[-\lambda,\lambda]^m$. Then for any $k \leq m/2$, with probability $1-c_0\exp(-c_1k)$, there is a permutation $\pi$ of $\{1,...,m\}$, such that for every $k/2 \leq i \leq k$,
$$
\frac{|\tau_{\pi(i)}|}{\lambda} \leq \frac{2i}{m}.
$$
\end{Lemma}

\proof Fix $i \leq  m/2$. Note that $\lambda^{-1} \tau_j$, $j=1,\ldots,m$, are independent and distributed uniformly in $[-1,1]$. Hence, for any $1\leq j\leq m$,
$$
\PP\left( \frac{|\tau_j|}{\lambda} \leq \frac{2i}{m}\right) = \frac{2i}{m}.
$$
By a binomial estimate, with probability at least $1-2\exp(-c_0 i)$,
$$
\left| \left\{j\in\{1,\ldots,m\} : \frac{|\tau_j|}{\lambda} \leq \frac{2i}{m} \right\} \right| \geq i.
$$
Now fix $k \leq m/2$, and the claim follows by the union bound over $k/2 \leq i \leq k$.
\endproof

Let us now complete the proof of Theorem~\ref{thm:main-B-2}. Let $T$ be a convex body. Fix $\delta>0$ and observe that for any $t \in T \cap \delta B_2^n$ we have that $f(t)=\sign(At + \tau)$ and $f(0)=\sign(\tau)$. Therefore,
$$
d_H(f(t),f(0)) = \sum_{i=1}^m \IND_{\{ \sign ( (At)_i + \tau_i) \not = \sign(\tau_i)\}}.
$$
Clearly, $\sign ( (A t)_i + \tau_i) \not = \sign(\tau_i)$ if and only if
$$
\tau_i \in [-(A t)_i,0) \cup [0,-(A t)_i).
$$

Let $k \leq m/2$ to be specified in what follows. By Lemma~\ref{lemma:lower-lambda}, with probability at least $1-2\exp(-c_1 k)$ there exists a permutation $\pi$ of $\{1,...,m\}$ (depending on the realization of $(\tau_i)_{i=1}^m$) for which
$$
\frac{|\tau_{\pi(i)}|}{\lambda} \leq \frac{2i}{m} \ \ \ {\rm for \ every } \ \ \frac{k}{2} \leq i \leq k.
$$
Conditioned on that event (which depends only on $\tau$), there is an absolute constant $c_2$ such that
$$
\left(\sum_{i=k/2}^{k} \tau_{\pi(i)}^2 \right)^{1/2} \leq c_2 \lambda \frac{k^{3/2}}{m}.
$$
Let $I=\{ \pi(i) : k/2 \leq i \leq k\}$ and let $P_I$ be the projection operator onto ${\rm span}(e_i)_{i \in I}$. Then $P_I A$ is a $\frac{k}{2} \times n$ standard gaussian matrix, and by Theorem \ref{thm:DvorMil}, with probability at least
$$
1-2\exp(-c_3d^*(T \cap \delta B_2^n)),
$$
we have that $P_I A (T \cap \delta B_2^n)$ contains the Euclidean ball
$$
c_4 \ell_*(T \cap \delta B_2^n) B_2^{|I|}
$$
as long as
\begin{equation}
\label{eqn:condonk1}
k \leq c_5 d^*(T \cap \delta B_2^n) \sim \frac{\ell_*^2(T \cap \delta B_2^n)}{\delta^2}.
\end{equation}
Hence, if in addition
\begin{equation}
\label{eqn:condonk2}
c_2 \lambda \frac{k^{3/2}}{m} \leq c_4 \ell_*(T \cap \delta B_2^n),
\end{equation}
then there is some $t^* \in T \cap \delta B_2^n$ such that for every $i \in I$
$$
\tau_i \in [-(A t^*)_i,0) \cup [0,-(A t^*)_i).
$$
In particular,
$$
d_H(f(t^*),f(0)) \geq \frac{k}{2} \ \ \  \text{and} \ \ \  \|t^*-0\|_2 \leq \delta.
$$
The required conditions \eqref{eqn:condonk1}, \eqref{eqn:condonk2}, and $k\leq m/2$ are all satisfied for
$$
k \sim \min\left\{ \left(\frac{m}{\lambda} \ell_*(T \cap \delta B_2^n)\right)^{2/3}, \frac{\ell_*^2(T \cap \delta B_2^n)}{\delta^2}\right\},
$$
thanks to the assumed lower bound on $\lambda$. 

Therefore,
$$
\frac{\lambda}{m} d_H(f(t^*),f(0)) \geq \frac{\lambda k}{2m} \geq c_6 \lambda \min\left\{ \frac{\ell_*^{2/3}(T \cap \delta B_2^n)}{\lambda^{2/3} m^{1/3}},\frac{\ell_*^2(T \cap \delta B_2^n)}{m \delta^2}\right\}.
$$
and if
$$
m \leq c_7\lambda \frac{\ell_*^2(T \cap \delta B_2^n)}{\delta^3}
$$
it follows that
$$
\left| \sqrt{2\pi} \frac{\lambda}{m} d_{H}(f(t^*),f(0)) - \|t^*-0\|_2 \right| \geq 2\delta.
$$
The claim is now evident by noticing that if $T$ is a convex body, then $T-T = 2T$.
\endproof

\subsection{The lower bound on $\lambda$}
The lower bound in Theorem \ref{thm:main-B-2} comes with a caveat: that
$$
\lambda \sqrt{m} \gtrsim \ell_*((T-T) \cap \delta B_2^n).
$$
Let us show that this caveat is not very restrictive.
\begin{Lemma} \label{lemma:lower-on-lambda}
There is an absolute constant $c$ and for every $\beta>0$ there is a constant  $\gamma=\gamma(\beta) \geq 1$ such that the following holds.
Assume that
\begin{equation} \label{eq:in-lemma-lower-on-lambda}
\ell_*( (T-T) \cap \lambda B_2^n) > \gamma \ell_*( (T-T) \cap \delta B_2^n).
\end{equation}
If, with probability at least $0.9$, the mapping $t \mapsto \sign(At+\tau)$ is a $\delta$-embedding of a convex body $T$, then
\begin{equation}
\label{eq:in-lemma-lower-on-lambda-conclusion}
\beta \lambda \sqrt{m} \geq \ell_*((T-T)\cap \delta B^n_2).
\end{equation}
\end{Lemma}
It follows that if the function $r \mapsto \ell_*( (T-T) \cap r B_2^n)$ exhibits a regular decay around $r=\delta$, then the condition in Lemma \ref{lemma:lower-on-lambda} is satisfied once $\lambda \gtrsim \delta$. As is explained in Appendix~\ref{sec:minShift}, the stronger condition $\lambda\gtrsim {\cal R}(T)$ is needed for the mapping $t \mapsto \sign(At+\tau)$ to have any chance of separating close points on rays $\{\alpha x : \alpha \geq 0\}$.

\begin{Remark}
The proof of Theorem \ref{thm:main-B-2} shows that the value of $\beta$ that is needed is just an absolute constant. Therefore, $\gamma$ in Lemma \ref{lemma:lower-on-lambda} may be taken to be an absolute constant as well.
\end{Remark}

\vskip0.3cm
The proof of Lemma \ref{lemma:lower-on-lambda} is based on the behaviour of the function
$$
\psi(r)=\frac{1}{\sqrt{m}r}\ell_*( (T-T) \cap r B_2^n).
$$
On the one hand, the map $r \mapsto \ell_*( (T-T) \cap r B_2^n)$ is increasing. On the other hand, it is standard to verify that $\psi$ is continuous and decreasing in $(0,\infty)$. Also, it tends to $0$ at infinity and to $\sqrt{n/m}$ at $0$ (because $T-T=2T$ is a convex body and in particular has a nonempty interior). Also note that by \eqref{eq:in-lemma-lower-on-lambda}, $\delta < \lambda/\gamma$. Indeed, \eqref{eq:in-lemma-lower-on-lambda} and the monotonicity of $r \mapsto \ell_*( (T-T) \cap r B_2^n)$ imply that $\delta \leq \lambda$ (as $\gamma\geq 1$). The monotonicity of $\psi$ therefore implies that $\psi(\lambda)\leq \psi(\delta)$ and hence
$$
\ell_*((T-T) \cap \lambda B_2^n) \leq \frac{\lambda}{\delta} \ell_*((T-T) \cap \delta B_2^n) < \frac{\lambda}{\delta} \cdot \frac{1}{\gamma} \ell_*((T-T) \cap \lambda B_2^n).
$$
The monotonicity of $\psi$ and the fact that $\delta<\lambda/\gamma$ are used frequently in the proof of Lemma~\ref{lemma:lower-on-lambda}.

\proof  Let $c_0$ be an absolute constant whose value is specified in what follows.  The key to the proof is to show that if $\rho$ satisfies
\begin{equation} \label{eq:lower-fixed-1}
\ell_*((T-T)\cap \delta B_2^n) \leq c_0 \frac{\rho}{\gamma} \sqrt{m},
\end{equation}
and
\begin{equation} \label{eq:lower-fixed-2}
\ell_*((T-T)\cap \rho B_2^n) \geq  c_0 \rho \sqrt{m}
\end{equation}
then $\lambda \geq c^\prime \rho$ for an absolute constant $c^\prime$. Once we establish that, the claim follows immediately: by \eqref{eq:lower-fixed-1}, we have that
$$
\ell_*((T-T)\cap \delta B^n_2) \leq c_0\frac{\rho}{\gamma} \sqrt{m} \leq  \frac{c_0}{c^\prime \gamma} \lambda\sqrt{m} \leq \beta \lambda \sqrt{m},
$$
provided that $\gamma \geq c_0/(c^\prime \beta)$.

The starting point of the proof is, once again, the Dvoretzky-Milman Theorem. By \eqref{eq:lower-fixed-2}, the Dvoretzky-Milman dimension of $(T-T)\cap \rho B_2^n$ satisfies
$$
d^*((T-T)\cap \rho B_2^n) \geq c_0^2 m.
$$
In particular, if $c_0$ is a sufficiently large absolute constant and $A:\R^n \to \R^m$ is the standard gaussian matrix, then with probability at least $0.9$,
$$
c_1 \ell_*((T-T)\cap \rho B_2^n) B_2^m \subset A((T-T)\cap \rho B_2^n)
$$
for an absolute constant $c_1$. For the rest of the proof, the values of $c_0$ and $c_1$ remain unchanged.

Before proceeding to the heart of the proof, let us show that we may assume without loss of generality that a suitable choice of $\rho$, satisfying \eqref{eq:lower-fixed-1} and \eqref{eq:lower-fixed-2}, exists.

Consider two cases: if $c_0 \geq \sqrt{n/m}$ then trivially, for any $r>0$, $\psi(r) \leq c_0$. In other words, for every $r>0$,
$$
\ell_*( (T-T) \cap r B_2^n) \leq c_0 r \sqrt{m}.
$$
In particular, for $r=\delta$,
$$
\ell_*( (T-T) \cap \delta B_2^n) \leq c_0 \delta \sqrt{m} \leq \beta \lambda \sqrt{m},
$$
when $\lambda \geq c_0\delta/\beta$. Recall that $\lambda > \gamma \delta$, and thus it suffices that $\gamma \geq c_0/\beta$ for the wanted estimate \eqref{eq:in-lemma-lower-on-lambda-conclusion} to hold.

Turning to the second case, when $c_0 \leq \sqrt{n/m}$, then since $\psi$ is continuous and decreasing from $\sqrt{n/m}$ to $0$, there is some $\rho^*$ for which
\begin{equation}
\label{eqn:rhoStarDef}
\ell_*((T-T)\cap \rho^* B_2^n) = c_0\rho^* \sqrt{m}.
\end{equation}
In particular, $\rho^*$ satisfies \eqref{eq:lower-fixed-2}.

Next, observe that either $\rho^*$ satisfies \eqref{eq:lower-fixed-1} and also $\lambda \leq \rho^*$, or, alternatively, the wanted estimate \eqref{eq:in-lemma-lower-on-lambda-conclusion} holds. Indeed, if $\lambda \geq \rho^*$ there is nothing to prove: in that case, using our assumption \eqref{eq:in-lemma-lower-on-lambda}, $\psi(\lambda)\leq \psi(\rho^*)$, and \eqref{eqn:rhoStarDef} we find
$$
\ell_*( (T-T) \cap \delta B_2^n) < \frac{1}{\gamma} \ell_*( (T-T) \cap \lambda B_2^n) \leq \frac{c_0}{\gamma} \lambda \sqrt{m} \leq \beta \lambda \sqrt{m}
$$
provided that $\gamma \geq c_0/\beta$. On the other hand, consider the case $\lambda \leq \rho^*$. If $\rho^*$ does not satisfy \eqref{eq:lower-fixed-1} then by \eqref{eqn:rhoStarDef} and \eqref{eq:in-lemma-lower-on-lambda},
$$
\ell_*((T-T)\cap \delta B_2^n) > c_0\frac{\rho^*}{\gamma} \sqrt{m} = \frac{1}{\gamma} \ell_*((T-T)\cap \rho^* B_2^n) \geq \frac{1}{\gamma} \ell_*((T-T)\cap \lambda B_2^n) > \ell_*((T-T)\cap \delta B_2^n),
$$
which is impossible.

Therefore, we may assume for the remainder of the proof that $c_0 \leq \sqrt{n/m}$, that $\rho^*$ satisfies both \eqref{eq:lower-fixed-1} and \eqref{eq:lower-fixed-2}, and that $\lambda \leq \rho^*$. To complete the proof we require three observations.

First, $\delta \leq \rho^*/\gamma$ since $\delta \leq \lambda \leq \rho^*$ and $\delta \leq \lambda/\gamma$.

Second, by applying Theorem~\ref{thm:gaussian-good-position} (with $k=m$) to the set
$$
(T-T) \cap \delta B_2^n=2T \cap \delta B_2^n,
$$
and using \eqref{eq:lower-fixed-1}, there is an absolute constant $c_2$ such that with probability at least $0.99$, for any $t \in T$,
$$
A\left( \left(t + \{z \in T : \|t-z\|_2 \leq \delta\} \right) \right) \subset At + c_2 \left(\delta+\frac{\rho^*}{\gamma}\right) \sqrt{m} B_2^m \subset At + c_3 \frac{\rho^*}{\gamma} \sqrt{m} B_2^m,
$$
where we have used that $\delta \leq \rho^*/\gamma$.

Third, as noted previously, by the choice of $c_0$ and the Dvoretzky-Milman Theorem, with probability at least $0.99$,
$c_1 \ell_*((T-T)\cap \rho^* B_2^n) B_2^m \subset A((T-T)\cap \rho^* B_2^n)$. Moreover, $T-T=2T$ and $2T \cap \rho^* B_2^n = 2(T \cap (\rho^*/2)B_2^n)$. Therefore,
$$
c_1\cdot c_0 \rho^* \sqrt{m} B_2^m = c_1 \ell_*((T-T)\cap \rho^* B_2^n) B_2^m \subset 2 A(T \cap (\rho^*/2) B_2^n).
$$
Now set $c_4=\frac{c_1}{2}c_0$ and note that $c_4$ is an absolute constant. If $\tau$ is supported in the interval $[-c_4\rho^*/10,c_4 \rho^*/10]$, then any two points in the same ``quadrant" of the set
$$
V=c_4 \rho^* \sqrt{m} B_2^{m} \cap \left\{ |x_i| \geq \frac{c_4 \rho^*}{10}, {\rm for \ every \ } 1 \leq i \leq m \right\},
$$
are indistinguishable for any realization of the function $y \mapsto \sign(y+\tau)$. Clearly, $V \subset A(T \cap (\rho^*/2) B_2^n)$ and the diameter of each of its quadrants is at least $(c_4/10) \rho^* \sqrt{m}$.

Fix a realization of $A$ in the intersection of the two events stated above---which has probability at least $0.98$. On that event, for every $t \in T$,
$$A\left( \left(t + \{z \in T : \|t-z\|_2 \leq \delta\} \right) \right)$$
has diameter at most $2c_3 \frac{\rho^*}{\gamma} \sqrt{m}$. On the other hand, any quadrant ${\cal Q}$ of $V$ has diameter at least $(c_4/10) \rho^* \sqrt{m}$. Hence, we can pick $y_1,y_2$ in ${\cal Q}$ with $\|y_1-y_2\|_2\geq (c_4/10) \rho^* \sqrt{m}$. Since $V \subset A(T \cap (\rho^*/2) B_2^n)$, $y_1=A t_1$ and $y_2=A t_2$ for some $t_1,t_2\in T$. If $\gamma$ is large enough so that 
$$2c_3 \frac{\rho^*}{\gamma}<(c_4/10) \rho^*$$  
then $y_1,y_2$ cannot be contained in $A\left(t + \{z \in T: \|t-z\|_2 \leq \delta\} \right)$ for a single $t \in T$; in particular, $\|t_1-t_2\|_2>\delta$. Finally, if $\lambda\leq c_4 \rho^*/10$, then since $y_1,y_2$ belong to the same quadrant of $V$ it would imply that $\sign(At_1+\tau)=\sign(At_2+\tau)$, violating the assumption that $t \mapsto \sign(At+\tau)$ is a $\delta$-embedding of $T$ with probability at least $0.9$. Thus, we must have $\lambda \geq \frac{c_4}{10}\rho^*$, and $c_4/10$ is an absolute constant---as required.    
\endproof

\subsection{A counterexample}
Next, let us explore the conjecture that $m \sim \ell_*^2(T)/\delta^2$ suffices to ensure that $f(t)=\sign(At+\tau)$ is, with high probability, a $\delta$-embedding of an arbitrary $T \subset \R^n$. We show that this conjecture is false, and that sometimes one in fact needs $m \geq c\ell_*^2(T)/\delta^3$.

Fix $0<\eps<\delta/2$, let $r \in \N$, set $I=\{1,...,r\}$ and put $J=\{r+1,...,\eta r/\eps^2\}$ for an absolute constant $\eta$ to be specified in what follows. Let $n=|I|+|J|$ and denote by $B_2^I$ and $B_2^J$ the Euclidean balls supported on ${\rm span}(e_i, i \in I)$ and ${\rm span}(e_j : j \in J)$ respectively.
Set
$$
{\cal E}=B_2^I \otimes \eps B_2^J = \left\{ (x,y): x \in B_2^I, \ y \in \eps B_2^J \right\} \subset \R^{n},
$$
and note that ${\cal E}$ is a convex body. 

Consider our random embedding $f:\R^n \to \{-1,1\}^m$, $f(t)=\sign(At+\tau)$. Observe that for any $0<\delta<\frac{1}{4}$, ${\cal E}$ contains two points $x,y$ located on a ray emanating from the origin such that $\|x\|_2=1$, $\|y\|_2<\|x\|_2$, and $\|x-y\|_2=2\delta$. As we show in detail in Appendix~\ref{sec:minShift}, this causes $f$ to fail to be a $\delta$-embedding of ${\cal E}$ with probability at least $0.9$ if $\lambda\leq c$, where $c$ is an absolute constant. Thus, to have any hope of obtaining a $\delta$-embedding, one needs to assume that $\lambda>c$, and $c$ is of the order of the diameter of ${\cal E}$. We will make a slightly stronger assumption: that $\lambda \geq C \sqrt{\log(e/\delta)}$; it corresponds to the context of Theorem~\ref{thm:main-A}.

Let us now assume that $m$ is such that, with probability at least $0.9$, $f$ is a $\delta$-embedding of ${\cal E}$. Using the notation of the previous section, the proof of Theorem \ref{thm:main-B-2} shows that for an absolute constant $\beta$, $\lambda$ must satisfy that
\begin{equation} \label{eq:in-counter-1}
\beta \lambda \sqrt{m} \geq \ell_*( ({\cal E}-{\cal E}) \cap \delta B_2^n).
\end{equation}
By Lemma \ref{lemma:lower-on-lambda}, it suffices to verify that
\begin{equation} \label{eq:in-counter-2}
\ell_*(({\cal E}-{\cal E}) \cap \lambda B_2^n) > \gamma(\beta) \ell_*(({\cal E}-{\cal E}) \cap \delta B_2^n)
\end{equation}
for a constant $\gamma=\gamma(\beta)$ (and thus, an absolute constant as well) to ensure that \eqref{eq:in-counter-1} holds.

Without loss of generality, we have that $\lambda \geq {\cal R}({\cal E})=\sup_{t \in {\cal E}} \|t\|_2=1$, and to verify \eqref{eq:in-counter-2} it is enough to ensure that
\begin{equation} \label{eq:in-counter-3}
\ell_*(2{\cal E}) \geq \gamma(\beta) \ell_*( 2{\cal E} \cap \delta B_2^n).
\end{equation}
Clearly, \eqref{eq:in-counter-3} holds if $\ell_*(B_2^I) \geq \gamma(\beta) (\ell_*(\delta B_2^I)+\ell_*(2\eps B_2^J))$, i.e., that 
$$
\sqrt{r} \geq c \gamma(\beta) (\delta\sqrt{ r}+\eps \sqrt{r}\sqrt{1+\eta/\eps^2}).
$$ 
It follows that if $\delta$ and $\eta$ are sufficiently small (depending on the absolute constant $\gamma(\beta)$), the wanted estimate holds.

Now we may assume that $\eta$ is a suitable, sufficiently small absolute constant and that $\delta$ is also small enough---ensuring that \eqref{eq:in-counter-3} holds. Set $\theta \sim \delta/\log(e\lambda/\delta)$, and by the choice of $\lambda$ we have that $\theta \sim \delta/\log(e/\delta)$.

\begin{Theorem} \label{thm:counter-example}
There are absolute constants $c_0$,$c_1$, and $c_2$ such that the following holds.
Let $\delta\leq c_0$ and let $A:\R^n \to \R^m$ be the standard gaussian matrix. If $m \leq c_1 \sqrt{\log(e/\delta)}\frac{\ell_*^2({\cal E})}{\delta^3}$,
then with probability at least $0.9$ there are $x,y \in {\cal E}$ such that
$$
\left| \sqrt{2\pi} \frac{\lambda}{m} d_{H}(f(x),f(y)) - \|x-y\|_2 \right| \geq 2\delta.
$$
Moreover, if $m \geq c_2 \sqrt{\log(e/\delta)}\frac{\ell_*^2({\cal E})}{\delta^3}$ then with probability at least 0.9
$$
\sup_{x,y \in {\cal E}} \left| \sqrt{2\pi} \frac{\lambda}{m} d_{H}(f(x),f(y)) - \|x-y\|_2 \right| \leq \delta.
$$
\end{Theorem}

\proof Set $\eps=\theta/2$ and let $U \subset B_2^I\otimes 0$ be a minimal $\theta/2$-cover of $B_2^I$ with respect to the Euclidean norm. Thus, $U$ is also a $\theta$-cover of ${\cal E}$: every $x \in  {\cal E}$ is supported on $I \cup J$, and there is some $u \in U$ such that
$$
\|x-u\|_2 = \|P_I x -u\|_2 + \|P_Jx\|_2 \leq \theta/2 + \eps \leq \theta.
$$
Moreover, by a standard volumetric estimate,
$$
\log {\cal N}({\cal E}, \theta) \leq \log {\cal N}(B_2^I, \theta/2) \lesssim r \log(5/\theta).
$$
Next, note that $\ell_*({\cal E})$, $\ell_*(B_2^I)$ and $\ell_*(\eps B_2^J)$ are all equivalent to $\sqrt{r}$. Moreover, using the choices of $\eps$ and $\theta$, it is straightforward to verify that there is an absolute constant $c_1$ such that 
$$
c_1\ell_*({\cal E}) \leq \ell_*(\eps B_2^J) \leq \ell_*({\cal E} \cap \theta B_2^n) \leq \ell_*({\cal E} \cap \delta B_2^n) \leq \ell_*({\cal E}).
$$
Hence, by Theorem \ref{thm:main-B-2}
$$
m \geq c_2 \sqrt{\log(e/\delta)} \frac{\ell_*^2({\cal E})}{\delta^3}
$$
is a necessary condition for ensuring that with high probability, $t \mapsto \sign(At+\tau)$ is a $\delta$-embedding of ${\cal E}$ in $\{-1,1\}^m$.

At the same time, recalling that $\lambda\sim \sqrt{\log(e/\delta)}$, we see that the second term in the upper bound \eqref{eqn:main-ABdm} from Theorem \ref{thm:main-A} (featuring the gaussian mean-width) is dominant for $\delta$ that is sufficiently small.  Thus, the lower bound is matched by the upper one from Theorem \ref{thm:main-A}.

\endproof

\begin{Remark}
(Possible logarithmic improvement) If $c<\lambda\leq C \sqrt{\log(e/\delta)}$, then Theorem~\ref{thm:main-A} is not applicable for $T={\cal E}$. However, Theorem \ref{thm:main-B-2} shows in this context that $f$ fails to be a $\delta$-embedding with probability at least $0.9$ if $m \leq c\frac{\ell_*^2({\cal E})}{\delta^3}$. We leave the question whether Theorem \ref{thm:main-A} continues to hold if $\lambda\geq C' {\cal R}(T)$ unanswered. If the answer is positive, then that would yield matching bounds for $T={\cal E}$ in the slightly improved setting $\lambda \sim {\cal R}({\cal E})$.
\end{Remark}

\appendix

\section{Proof of Lemma \ref{lem:exp}} \label{app:proof-lem:exp}
We will use the following simple observation. 
\begin{Lemma}
\label{lem:sepProbTau}
Let $a,b\in \R$. If $\tau$ is uniformly distributed on $[-\lambda,\lambda]$, then 
\begin{align*}
2\lambda\bP(\sign(a+\tau)\neq \sign(b+\tau)) & = |a-b| 1_{\{|a|\leq \lambda, |b|\leq \lambda\}} + 2\lambda 1_{\{a>\lambda, b<-\lambda\}} + 2\lambda 1_{\{a<-\lambda, b>\lambda\}}  \\
& \qquad + (\lambda-a) 1_{\{b>\lambda,|a|\leq \lambda\}} + (\lambda-b) 1_{\{a>\lambda,|b|\leq \lambda\}} \\ 
& \qquad + (\lambda+a) 1_{\{b<-\lambda,|a|\leq \lambda\}} + (\lambda+b) 1_{\{a<-\lambda,|b|\leq \lambda\}}. 
\end{align*}    
\end{Lemma}
\proof
The result follows easily by distinguishing cases. Write 
$$\alpha:= 2\lambda\bP(\sign(a+\tau)\neq \sign(b+\tau)).$$ 
Suppose first that $|a|,|b|\leq \lambda$. If $a<b$, then 
\begin{equation*}
\PP(\sign (a+\tau)\neq \sign (b+\tau))=\PP(a<-\tau \leq b)=\frac{b-a}{2\lambda}.
\end{equation*}
By reversing the roles of $a$ and $b$ we find $\alpha= |b-a|$.\par
Next, if either $a,b>\lambda$ or $a,b<-\lambda$, then $\sign(a+\tau)=\sign(b+\tau)$ and hence $\alpha=0$.\par
If $a>\lambda$ and $b<-\lambda$, then $\sign(a+\tau)=1\neq -1 =\sign(b+\tau)$ and so $\alpha=2\lambda$. Reversing the roles of $a$ and $b$, we see that the same holds if $b>\lambda$ and $a<-\lambda$.\par
If $a>\lambda$ and $|b|\leq \lambda$, then $\sign(a+\tau)=1$ and so
$$\alpha = 2\lambda \bP(\sign(b+\tau)=-1) = 2\lambda\bP(b+\tau<0) = 2\lambda\frac{-b-(-\lambda)}{2\lambda} = \lambda-b.$$
Reversing the roles of $a$ and $b$, we see that $\alpha=\lambda-a$ if $b>\lambda$ and $|a|\leq \lambda$.\par
Finally, suppose that $a<-\lambda$ and $|b|\leq \lambda$. Then $\sign(a+\tau)=-1$ and hence
$$\alpha = 2\lambda \bP(\sign(b+\tau)=1) = 2\lambda\bP(b+\tau\geq 0) = 2\lambda\frac{\lambda-(-b)}{2\lambda} = \lambda+b.$$
Reversing the roles of $a$ and $b$ shows that $\alpha=\lambda+a$ if $b<-\lambda$ and $|a|\leq \lambda$.
\endproof
We will now prove Lemma~\ref{lem:exp} using Lemma~\ref{lem:sepProbTau}. We assume without loss of generality that $a<b$ and distinguish cases. If $|a|\leq \lambda$ and $|b|\leq \lambda$ then clearly $\phi_\lambda(a)+\phi_\lambda(b)=0$ and
\begin{equation*}
2\lambda \PP(\sign (a+\tau) \neq \sign (b+\tau))=|a-b|.
\end{equation*}
Next, consider the case in which either $|a|>\lambda$ or $|b|>\lambda$ (or both), and we distinguish
cases according to the value of $a$:
\begin{itemize}
\item [1)] \underline{$a<-\lambda$.} There are three different ranges of $b$:
\\
If $b<-\lambda$ then 
\begin{align*}
& \left|2\lambda \PP(\sign(a+\tau)\neq \sign(b+\tau))- |a-b| \right|
\\
& \qquad  =  |a-b| \leq |a+\lambda|+|b+\lambda| =-(a+\lambda)-(b+\lambda) =\phi_\lambda(a)+\phi_\lambda(b).
\end{align*}
If $|b|\le \lambda$ then 
\begin{align*}
& \left| 2\lambda \PP(\sign(a+\tau) \neq \sign(b+\tau))- |a-b| \right|
\\
& \qquad = | b+\lambda -(b-a) | = \phi_\lambda(a) = \phi_\lambda(a)+\phi_\lambda(b).
\end{align*}
Finally, if $\lambda<b$ then 
\begin{align*}
& \left|2\lambda \PP(\sign(a+\tau) \neq \sign(b+\tau))- |a-b|\right|
\\
& \qquad = |2\lambda  - (b-a)| =(b-a)-2\lambda  =-a-\lambda + b-\lambda  =\phi_\lambda(a)+\phi_\lambda(b).
\end{align*}
\item [2)] \underline{$|a|\leq \lambda$.} In this case, $\lambda<b$ and hence
\begin{align*}
& \left|2\lambda \PP(\sign(a+\tau) \neq \sign (b+\tau))- |a-b| \right|
\\
& \qquad = | \lambda-a  - (b-a)  | =b-\lambda =|b|-\lambda =\phi_\lambda(b) = \phi_\lambda(a)+\phi_\lambda(b).
\end{align*}
\item [3)] \underline{$\lambda <a$.} Then $\lambda <b$ and so
\begin{align*}
& \left|2\lambda \PP(\sign(a+\tau)\neq \sign(b+\tau))- |a-b| \right|
\\
& \qquad = b-a = b-\lambda - (a-\lambda) = \phi_\lambda(b) - \phi_\lambda(a) \leq \phi_\lambda(a)+\phi_\lambda(b).
\end{align*}
\end{itemize}
\endproof

\section{Minimal shift}
\label{sec:minShift}

As before, let $f:\R^n\to \{-1,1\}^m$ be the map $f(t)=\sign(At+\tau)$, where $A:\R^n\to\R^m$ is the gaussian matrix, $\tau$ is uniformly distributed in $[-\lambda,\lambda]^m$ and $A$ and $\tau$ are independent. In this section we demonstrate that, in general, $f$ can easily fail to be a $\delta$-embedding of a set $T$ unless the shift parameter $\lambda$ is of the order of $\mathcal{R}(T)$. To see this, let us first note that $f$ cannot be a $\delta$-embedding unless $m$ exceeds a minimal threshold. Indeed, suppose that $f$ is a $\delta$-embedding of a two-point set $\{x,y\}$ and that $\|x-y\|=2\delta$. Then, clearly 
$$\frac{\sqrt{2\pi}\lambda}{m}d_H(f(x),f(y))\geq\delta \qquad \text{and} \qquad \frac{\sqrt{2\pi}\lambda}{m}d_H(f(x),f(y))\leq 3\delta.$$
The first inequality implies $d_H(f(x),f(y))\geq 1$ and hence $m\geq \lambda\sqrt{2\pi}/(3\delta)$ by the second inequality. We will therefore assume $m\geq c\lambda/\delta$ in what follows.\par 
The following Lemma formalizes the key `difficult' geometric arrangement that requires large shifts.   
\begin{Lemma}
\label{lem:badTwoPointSet}
There exist absolute constants $c_1,c_2>0$ such that the following holds. Let $\delta<\frac{1}{4}$ and let $x,y$ be two points located on a ray emanating from the origin such that $\|x\|_2\geq 4\delta$, $\|y\|_2<\|x\|_2$ and $\|x-y\|_2=2\delta$. If $\lambda\leq c_1\|x\|_2$ and $m\geq c_2\lambda/\delta$, then $f$ is not a $\delta$-embedding of $\{x,y\}$ with probability at least $0.9$.
\end{Lemma}
\proof As we noted before, if $f$ is a $\delta$-embedding of $\{x,y\}$, then  
\begin{equation}
\label{eqn:deltaEmbedCons}
\frac{\sqrt{2\pi}\lambda}{m}d_H(f(x),f(y))\geq\delta
\end{equation} 
Set 
$$p=\bP(\sign(\langle g,x\rangle+\tau)\neq \sign(\langle g,y\rangle+\tau)),$$
then by Chernoff's inequality, for any $u>mp$,
$$\bP(d_H(f(x),f(y))\geq u)\leq e^{-mp} \left(\frac{emp}{u}\right)^u.$$ 
We will now show that $\sqrt{2\pi}\lambda p\leq \delta/e^2$ under the stated conditions. We can then take $u=\delta m/(\lambda\sqrt{2\pi})$ to see that \eqref{eqn:deltaEmbedCons} is violated with probability at least $0.9$ if $m\geq c_2\lambda/\delta$.\par
By applying Lemma~\ref{lem:sepProbTau} with $a=\langle g,x\rangle$ and $b=\langle g,y\rangle$, we find
\begin{align}
\label{eqn:sepProbTauGaussian}
& 2\lambda\bP_{\tau}(\sign(\langle g,x\rangle+\tau)\neq \sign(\langle g,y\rangle+\tau)) \nonumber\\
& \qquad = |\langle g,x-y\rangle| 1_{\{|\langle g,x\rangle|\leq \lambda, |\langle g,y\rangle|\leq \lambda\}} + 2\lambda (1_{\{\langle g,x\rangle>\lambda, \langle g,y\rangle<-\lambda\}} + 1_{\{\langle g,x\rangle<-\lambda, \langle g,y\rangle>\lambda\}})  \nonumber\\
& \qquad \qquad + (\lambda-\langle g,x\rangle) 1_{\{\langle g,y\rangle>\lambda,|\langle g,x\rangle|\leq \lambda\}} + (\lambda-\langle g,y\rangle) 1_{\{\langle g,x\rangle>\lambda,|\langle g,y\rangle|\leq \lambda\}} \nonumber \\ 
& \qquad \qquad + (\lambda+\langle g,x\rangle) 1_{\{\langle g,y\rangle<-\lambda,|\langle g,x\rangle|\leq \lambda\}}  + (\lambda+\langle g,y\rangle) 1_{\{\langle g,x\rangle<-\lambda,|\langle g,y\rangle|\leq \lambda\}}. 
\end{align}    
Since $x$ and $y$ are located on a ray emanating from the origin, 
$$1_{\{\langle g,x\rangle>\lambda, \langle g,y\rangle<-\lambda\}} + 1_{\{\langle g,x\rangle<-\lambda, \langle g,y\rangle>\lambda\}}=0.$$
Combining this with trivial estimates for the other terms in \eqref{eqn:sepProbTauGaussian}, we find
\begin{align*}
& 2\lambda\bP_{\tau}(\sign(\langle g,x\rangle+\tau)\neq \sign(\langle g,y\rangle+\tau)) \\
& \qquad \leq 3|\langle g,x-y\rangle| 1_{\{|\langle g,x\rangle|\leq \lambda\}} + 2|\langle g,x-y\rangle| 1_{\{|\langle g,y\rangle|\leq \lambda\}}
\end{align*}
Let $\gamma$ be a standard gaussian random variable, then $\langle g,v\rangle$ is identically distributed with $\gamma\|v\|_2$, for any $v\in \R^n$. Taking expectations with respect to $g$ and using the Cauchy-Schwarz inequality, we find 
\begin{align*}
& \sqrt{2\pi}\lambda\bP(\sign(\langle g,x\rangle+\tau)\neq \sign(\langle g,y\rangle+\tau)) \\
& \qquad \leq C \ \E\big(|\langle g,x-y\rangle| 1_{\{|\langle g,x\rangle|\leq \lambda\}} + |\langle g,x-y\rangle| 1_{\{|\langle g,y\rangle|\leq \lambda\}}\big) \\
& \qquad = 2C\delta ((\bP(|\langle g,x\rangle|\leq \lambda))^{1/2}+(\bP(|\langle g,y\rangle|\leq \lambda))^{1/2}).
\end{align*}
Note that $y=(1-\frac{2\delta}{\|x\|_2})x$ and hence 
$$\bP(|\langle g,y\rangle|\leq \lambda) = \bP(|\gamma|\leq \lambda/\|y\|_2) \leq \frac{\lambda}{\|y\|_2} = \frac{\lambda}{\|x\|_2-2\delta} \leq \frac{2\lambda}{\|x\|_2}$$ 
as $\|x\|_2\geq 4\delta$. Similarly,
$$\bP(|\langle g,x\rangle|\leq \lambda)\leq \frac{\lambda}{\|x\|_2}.$$
Thus, there is a constant $c_0>0$ such that if $\lambda\leq c_0 \|x\|_2$, then 
$$(\bP(|\langle g,x\rangle|\leq \lambda))^{1/2}+(\bP(|\langle g,y\rangle|\leq \lambda))^{1/2}\leq \frac{1}{2e^2C}.$$
\endproof

\section{Proof of Theorem~\ref{thm:gaussian-good-position}}
\label{sec:knormGaussian}

Let us write $\alpha\leq \beta$ if $\alpha \leq C\beta$ for an absolute constant $C>0$. Let $\|\cdot\|_{\psi_2}$ denote the subgaussian norm. Let $a_1^T,\ldots,a_m^T$ be the rows of $A$, let $\Sigma_{k,m}$ be the set of all $k$-sparse vectors in the Euclidean unit ball, and set $T':=\mathcal{R}(T)^{-1} T$. With this notation, we can write
$$\sup_{x\in T} \|Ax\|_{[k]} = \mathcal{R}(T)\sup_{(x,y)\in T'\times \Sigma_{k,m}} \langle Ax,y\rangle = \mathcal{R}(T)\sup_{(x,y)\in T'\times \Sigma_{k,m}} \sum_{i=1}^m \langle a_i,x\rangle y_i.$$
Note that the right hand side is the supremum of a stochastic process with mean zero and subgaussian increments. Indeed, for $(x,y), (\tilde{x},\tilde{y}) \in T'\times \Sigma_{k,m}$, Hoeffding's inequality for subgaussian random variables (see, e.g., \cite[Section 2.6]{Ver18}) yields
\begin{align*}
& \left\|\sum_{i=1}^m \langle a_i,x\rangle y_i - \sum_{i=1}^m \langle a_i,\tilde{x}\rangle \tilde{y}_i\right\|_{\psi_2} \\
& \qquad  =  \left\|\sum_{i=1}^m \langle a_i,x-\tilde{x}\rangle y_i + \sum_{i=1}^m \langle a_i,\tilde{x}\rangle (y-\tilde{y})_i\right\|_{\psi_2} \\
& \qquad \leq \left\|\sum_{i=1}^m \langle a_i,x-\tilde{x}\rangle y_i\right\|_{\psi_2} + \left\|\sum_{i=1}^m \langle a_i,\tilde{x}\rangle (y-\tilde{y})_i\right\|_{\psi_2} \\
& \qquad \lesssim \left(\sum_{i=1}^m \|\langle a_i,x-\tilde{x}\rangle y_i\|_{\psi_2}^2\right)^{1/2} + \left(\sum_{i=1}^m \|\langle a_i,\tilde{x}\rangle (y-\tilde{y})_i\|_{\psi_2}^2\right)^{1/2} \\
& \qquad \lesssim \left(\sum_{i=1}^m \|x-\tilde{x}\|_2^2 y_i^2\right)^{1/2} +  \left(\sum_{i=1}^m \|\tilde{x}\|_2^2 (y-\tilde{y})_i^2\right)^{1/2} \\
& \qquad \leq \|x-\tilde{x}\|_2 + \|y-\tilde{y}\|_2 \lesssim  \|(x,y)-(\tilde{x},\tilde{y})\|_2.  
\end{align*}
By \cite[Theorem 3.2]{Dir15} and Talagrand's majorizing measures theorem (see, e.g., \cite{Tal14}), for any $v\geq 1$ we find that with probability at least $1-e^{-v^2}$,
$$\sup_{x\in T} \|Ax\|_{[k]}\lesssim \mathcal{R}(T) (\ell_*(T'\times \Sigma_{k,m})+v) \lesssim  \ell_*(T) + \mathcal{R}(T)\sqrt{k\log(em/k)}+ \mathcal{R}(T)v.$$
Setting $v=u\sqrt{k\log(em/k)}$ yields the result.

\end{document}